\documentclass[11pt]{article}
\usepackage{a4}
\usepackage{amsmath}

\usepackage{amsfonts}
\usepackage{amssymb}
\usepackage{euscript}
\newenvironment{proof}{\noindent {\it Proof.~~}\ }{\
  \rule{1mm}{2mm}\medskip}
\newenvironment{proofof}[2]{\noindent {\it Proof of #1}~#2: \
}{~\rule{1mm}{2mm}\medskip}

\renewenvironment{itemize}{\begin{list}{$\bullet$}{\topsep=0mm\itemsep=0mm\labelwidth=10mm}}{\end{list}}

\newtheorem{theorem}{Theorem}
\newtheorem{lemma}[theorem]{Lemma}
\newtheorem{corollary}[theorem]{Corollary}

\title{On some subgroup chains related to Kneser's theorem}

\author{Yahya O. Hamidoune\thanks{Universit\'e Pierre et Marie Curie, Paris 6,
Combinatoire et Optimisation - case 189,
4 place Jussieu,
75252 Paris Cedex 05. {\tt yha@ccr.jussieu.fr} }
\and Oriol Serra\thanks{Universitat Polit\`ecnica de Catalunya,
Matem\`atica Aplicada IV,
Campus Nord - Edif. C3,
C. Jordi Girona, 1-3,
08034 Barcelona, Spain. {\tt oserra@mat.upc.es} } \and Gilles
Z\'emor\thanks{Universit\'e de Bordeaux 1, Institut de Math\'ematiques
de Bordeaux, 351 cours de la Lib\'eration, 33405 Talence.
{\tt
zemor@math.u-bordeaux1.fr} }}

\date{March 28, 2007}

\begin{document}
\maketitle

\begin{abstract}

A recent result of Balandraud shows that for every subset $S$ of an
abelian group $G$ there
exists a non trivial subgroup $H$ such that
 $|TS|\le |T|+|S|-2$ holds only if  $H\subset Stab (TS)$. Notice
 that Kneser's Theorem only gives $\{0\}\neq Stab (TS)$.

 This strong form of Kneser's theorem follows from  some nice properties
 of a certain poset investigated by Balandraud.
 We consider an analogous poset for nonabelian groups and, by
 using classical tools from Additive Number Theory, extend some
 of the above results. In particular we obtain short proofs of 
 Balandraud's results in the abelian case.
\end{abstract}

\section{Introduction}

In order to avoid switching from multiplicative to additive notation,
all groups will be written multiplicatively.

Kneser's addition theorem states that if $S, T$ are finite subsets
of an abelian group $G$ then $|ST|\le |S|+|T|-2$ holds only if $ST$
is periodic (i.e, there is a non trivial subgroup $H$ such that
$HST=ST$.) Kneser's Theorem is a fundamental tool in Additive number
Theory. Proofs of this result may be found in
\cite{KEM,KN3,MAN,Nat,tv}.

 In all previously known
proofs of Kneser's Theorem, the subgroup $H$ depends crucially on both sets $S$
and $T$. With the goal of breaking this double dependence in $S$ and
$T$, Balandraud investigated in recent work \cite{bal,bal1}
the properties of a combinatorial poset that we now present.

Let $S$ be a finite subset containing $1$ of a group $G$. Following
Balandraud, let us define a  {\em cell} of $S$ as a finite subset $X$ such
that, for all $z\notin X$, it holds that $zS\not\subset XS$. This
notion is defined in \cite{bal,bal1} and it is equivalent to the
notion of nonextendible subset used in \cite{dav}. Throughout the
paper, by a cell we always mean a cell of $S$.

A cell $X$ is called a $u$-{\em cell} if $|XS|-|X|=u$.
 A $u$-cell with minimal cardinality is called  a $u$-{\em kernel}
 (of $S$).

 Balandraud showed that, for a finite set $S$ in an  abelian
group $G$, in the poset of $j$--cells  containing the unity ordered
by inclusion with $1\le j\le |S|-2$, the set of kernels form a chain
of subgroups. Moreover, if there exists a $u$--cell, then there is a
unique $u$--kernel containing the unit element which is contained in all
$u$--cells containing the unit element.

One of the consequences of this work  is a new proof and the
following strengthening of Kneser's Theorem:

\begin{theorem}[Balandraud]\label{thm:bal}
  For any non-empty finite subset $S$ of an abelian group $G$, 
  there exists a finite subgroup $H$ of $G$ such that for any
  finite subset $T$ of $G$ one of the following conditions hold~:
  \begin{itemize}
  \item $|TS|\geq |T|+|S|-1$
  \item $HTS=TS$ and $|TS|\leq |HS|+|HT|-|H|$
  \end{itemize}
\end{theorem}

As far as the authors are aware this is a surprising and strong
formulation that was not observed before and does not follow
straightforwardly from the classical forms of Kneser's Theorem.

The purpose of the present note is to give a short proof for
the nonabelian case that, in the poset
of $j$--cells  {\em that are subgroups} ordered by inclusion with
$0\le j\le |S|-1$, the set of kernels form a chain of subgroups.
Moreover, each $u$-kernel of this poset is unique and contained in
all $u$--cells of this poset.

{F}rom this statement  Kneser's theorem allows one to deduce
Balandraud's results for the abelian case, and in particular Theorem
\ref{thm:bal}. Kneser's Theorem has several equivalent forms.
We use the following one; see e.g \cite{KEM,Nat}:

\begin{theorem}[Kneser \cite{KN3}]
Let $G$ be an abelian group and $X, Y \subset G$ be finite subsets
such that $|XY|\le |X|+|Y|-2$. Then
$$|XY|=|HX|+|HY|-|H|,$$  where $H=stab (XY)=\{x : xXY=XY\}.$
\label{Kneser}
\end{theorem}

Our main tool is the following Theorem of
Olson\cite[Theorem 2]{ols}. We give an equivalent
formulation here where we use left--cosets instead of right--cosets.

\begin{theorem}[Olson \cite{ols}]
Let $X, Y$ be finite subsets of a group $G$, and let $H$ and $K$ be
subgroups such that  $HX=X$, $KY=Y$  and $KX\neq X$, $HY\neq Y$.
Then
$$|X\setminus  Y|+ |Y\setminus  X|\geq |H|+|K|-2|H\cap K|.$$
In particular either     $|X\setminus Y|\ge |H|-|H\cap K|$ or
$|Y\setminus X|\ge |K|-|H\cap K|$.
 \label{olson}
\end{theorem}

We shall use the following lemma.

\begin{lemma}[\cite{bal,bal1}]\label{lem:cell}
Let $G$ be a group and $1\in S \subset G$ be a finite subset. Then
the intersection of two cells $M_1, M_2$ of $S$ is a cell of $S$.
\end{lemma}

\begin{proof}
Let $x\notin M_1\cap M_2$. There is $i$ with $x\notin M_i$. Then
$xS\not\subset M_iS.$  Hence $xS \not\subset (M_1\cap M_2)S$.
\end{proof}

We can now state our main result, namely Theorem \ref{thm:mn} below.

\section{An application of Olson's Theorem}

Balandraud \cite{bal,bal1} proved  that, in the abelian case, the
set of kernels containing the unit element and ordered by inclusion is
a chain of subgroups.
In the non abelian case we can prove only that the set of kernels
that are {\em subgroups} forms a chain. The abelian case can then be
easily recovered,
 since Kneser's Theorem implies (as we shall see below)  
that a kernel containing the unit element is a
subgroup.

\begin{theorem}\label{thm:mn}
Let $S$ be a finite subset containing $1$ of a group $G$. Let $M$ be
a $u$--kernel of $S$ which is a subgroup. Let $N$ be a subgroup which is a
$v$--cell and suppose $u,v\le |S|-1$.
\begin{itemize}
\item[(i)]
If either $N$ is a $v$--kernel or $u=v$ then $M\subset N$ or
$N\subset M$.

\item[(ii)]If $N$ is a $v$--kernel and  $v\le u$ then $M\subset N$.
\end{itemize}
\end{theorem}

\begin{proof}
Suppose  that $M\not\subset N$ and $N\not\subset M$.  Note that,
since $M$ is a cell, if $NMS=MS$ then $NM= M$, thus $N\subset M$
against our assumption. Hence we may assume $NMS\neq MS$ and
similarly $MNS\neq NS$. By Theorem \ref{olson} we have one of the
two following cases.

{\bf Case } 1: $ |MS|-|(MS)\cap (NS)|=|(MS)\setminus (NS)|\ge
|M|-|M\cap N|.$ It follows that $|(M\cap N)S|-|M\cap N|\leq
|(MS)\cap (NS)|-|M\cap N|\leq|MS|-|M|.$ On the other hand we have
$u=|MS|-|M|<|S|\le |(M\cap N)S|.$ Since $|MS|-|M|$ is a multiple of
$|M\cap N|$ we  have
  $$u=|MS|-|M|=|(M\cap N)S|- |M\cap N|.$$
By Lemma \ref{lem:cell},
$M\cap N$ is a cell. Since $M$ is a  $u$--kernel, we have
$M\cap N=M,$ a contradiction.

{\bf Case } 2: $ |NS|-|(NS)\cap (MS)|=|(NS)\setminus (MS)|\ge
|N|-|N\cap M|.$ It follows that $|(N\cap M)S|-|N\cap M|\leq
|(NS)\cap (MS)|-|N\cap M|\leq|NS|-|N|.$ On the other hand we have
$|NS|-|N|<|S|\le |(N\cap M)S|.$ Since $|NS|-|N|$ is a multiple of
$|N\cap M|$ we have
\begin{equation}\label{eq:o2}
|NS|-|N|=|(N\cap M)S|- |N\cap M|.
\end{equation}

Assume first $u=v$. Then $u=|MS|-|M|=|NS|-|N|=|(N\cap M)S|- |N\cap
M|$. Since $M$ is a  $u$--kernel, we have $M\cap N=M,$ a
contradiction.

Assume that $N$ is a $v$--kernel. Then (\ref{eq:o2}) implies $N\cap
M=N,$ a contradiction. This proves $(i)$.

Assume now that $v\le u$. Suppose $M\not\subset N$. By $(i)$ we have
$N\subset M$, which implies in particular that
$|MS|-|M|$ is a multiple of $N$. Therefore,
from $u=|MS|-|M|<|S|\le |NS|$ we have
$u=|MS|-|M|\le |NS|-|N|=v$ which gives $u=v$. But then $M\not\subset N$
and $N\subset M$ imply $|N|<|M|$, and since $N$ is now a $u$--cell, this
contradicts $M$ being a $u$--kernel.
\end{proof}

We can now deduce Balandraud's description for kernels and cells :

\begin{corollary}[Balandraud \cite{bal,bal1}]
Let $G$ be an abelian group and $S \subset G$ be a finite subset
with $1\in S$. Let $M$ be a $u$--kernel of $S$ containing $1$ with
$1\le u\le |S|-2$. Then,
\begin{itemize}
\item[{\rm (i)}]
$M$ is a subgroup.
\item[{\rm (ii)}] Each  $u$-cell  is $M$--periodic.
\item[{\rm (iii)}] Each $v$--kernel with $u<v\leq |S|-2$
 is a proper subgroup of $M$.
\end{itemize}
\label{balthm}
\end{corollary}

\begin{proof} Let $X$ be a $u$-cell with $u\le |S|-2$.
By Kneser's Theorem, the inequality $|XS|-|X|=u\le |S|-2$ implies
\begin{equation}\label{eq:kn}
u=|XS|-|HX|=|HS|-|H|, \end{equation}
 where $H$ is the stabilizer of
$XS$. Since $X$ is a cell and $HXS=XS$, we have $X=HX$. Note that,
since $G$ is abelian,
$(\{y\}\cup H)S=HS$ implies $y\in Stab(HS)\subset Stab(XS)$, so that
$y\in H$. This observation and (\ref{eq:kn}) imply   that $H$ is an
$u$--cell. In particular, by taking $X=M$, the period $K$ of $MS$ is
a $u$--cell. Since $KMS=MS$ and $M$ is a $u$--cell, we have
$K\subset KM\subset M$. Since $M$ is a $u$--kernel we have $M=K$.
This proves (i).

Now let $H$ be the stabilizer of $XS$, where $X$ is a $u$--cell. As
shown in the preceding paragraph $H$ is also a $u$--cell. By
Theorem \ref{thm:mn} we have $M\subset H$ and thus $MH=H$. Since $X$
is a cell and $HXS=XS$, we have $X=HX=MHX$. Hence $X\subset MX\subset MHX =X$
implies $X=MX$. This proves (ii).

Finally, by (i), a $v$--kernel $N$ is a subgroup. By Theorem
\ref{thm:mn} we have $N\subset M$.
\end{proof}

{F}rom Corollary \ref{balthm}, one can  deduce Theorem \ref{thm:bal}.

\begin{proofof}{Theorem}{\ref{thm:bal}}
We may assume without loss of generality that $1\in S$.

{\bf Case 1:}
There is no $m$--cell for any $1\leq m\leq |S|-2$.
\begin{itemize}
\item either we have $|TS|\geq |S|+|T|-1$ for any non-empty finite
      $T$, in which case the theorem clearly holds with $H=\{1\}$.
\item or there exists some non-empty finite $T$ such that
      $|TS|\leq |S|+|T|-2$. Without loss of generality, we may also
      suppose $1\in T$. Now $T$ must be contained in an $m$--cell
      with $m\leq |S|-2$, but since no such cell exists for $1\leq m$,
      we have that $T$ itself must be a cell (a $0$-cell) i.e.
      $|TS|=|T|$. We therefore have $HT=TH=T=TS=HTS$ where
      $H$ is the (necessarily finite) subgroup generated by $S$.
      We have just proved that the theorem holds
      in this case with $H=\langle S\rangle$.
\end{itemize}

{\bf Case 2:}
There exists an $m$--cell with $1\leq m\leq |S|-2$.
We may therefore consider the largest integer $u\le |S|-2$
for which $S$ admits a $u$--cell. Let $H$ be the $u$--kernel
containing $1$. Note that $u\leq |S|-2$ implies that $H$ is different
from $\{1\}$. 
Now let $T$ be any finite non-empty subset such that $|TS|-|T|\leq
|S|-2$.
We shall prove that $HTS = TS$.

By adding elements to $T$ as long
as necessary, we can find a cell $X$ that contains $T$ and such that
$XS=TS$. Note that we then have $|XS|-|X|\leq |TS|-|T|\leq |S|-2$,
so that $X$ is a $v$--cell for some $v\leq u$.
By Corollary \ref{balthm} (ii) we have
$TS=XS=MXS=MTS$ where $M$ is the $v$-kernel containing $1$. By part (i)
of Corollary \ref{balthm}, $H$ is a subgroup of $M$ so that
$TS=XS=HTS$ as well.

Finally, $|ST|\leq |HS|+|HT|-|H|$ follows from $|ST|$ being a multiple
of $|H|$.
\end{proofof}

\end{document}